\documentclass[11pt]{article}
\usepackage[margin=1in]{geometry}
\usepackage{amsfonts,amsmath,amssymb}
\usepackage{algpseudocode}

\PassOptionsToPackage{obeyspaces}{url}
\usepackage[colorlinks=true,citecolor=blue,urlcolor=blue,linkcolor=blue,bookmarksopen=true]{hyperref}
\usepackage{amsthm}
\usepackage{enumitem}

\usepackage{tikz}

\usepackage{etoolbox}
\patchcmd{\thebibliography}{\leftmargin\labelwidth}{\leftmargin\labelwidth\addtolength\itemsep{-0.1\baselineskip}}{}{}

\author{Christopher Cox\thanks{Department of Mathematics, Iowa State University, Ames, IA, USA. \texttt{cocox@iastate.edu}. Supported in part through NSF RTG Grant DMS-1839918.}
\and Ryan R.\ Martin\thanks{Department of Mathematics, Iowa State University, Ames, IA, USA. \texttt{rymartin@iastate.edu}. Supported in part through Simons Collaboration Grants \#353292 and \#709641.}
\and Daniel McGinnis\thanks{Department of Mathematics, Iowa State University, Ames, IA, USA. \texttt{dam1@iastate.edu}. Supported in part through NSF RTG Grant DMS-1839918.}}

\title{Accumulation points of the edit distance function}
\date{}

\usepackage[nameinlink,sort]{cleveref}

\newtheorem{theorem}{Theorem}
\newtheorem{lemma}[theorem]{Lemma}

\newtheorem{conj}[theorem]{Conjecture}
\crefname{conj}{conjecture}{conjectures}

\newtheorem{claim}[theorem]{Claim}
\crefname{claim}{claim}{claims}

\newtheorem{prop}[theorem]{Proposition}
\crefname{prop}{proposition}{propositions}

\theoremstyle{definition}

\newtheorem{defn}[theorem]{Definition}
\crefname{defn}{definition}{definitions}

\newtheorem{remark}[theorem]{Remark}
\crefname{remark}{remark}{remarks}

\newtheorem{question}[theorem]{Question}
\crefname{question}{question}{questions}

\crefname{enumi}{part}{parts}

\newlist{proplist}{enumerate}{1}
\setlist[proplist,1]{label=\Roman*., ref=\Roman*}

\crefname{proplisti}{property}{properties}

\numberwithin{theorem}{section}

\makeatletter
\DeclareRobustCommand{\crefnosort}[1]{%
    \begingroup\@cref@sortfalse\cref{#1}\endgroup
}
\DeclareRobustCommand{\Crefnosort}[1]{%
    \begingroup\@cref@sortfalse\Cref{#1}\endgroup
}
\makeatother

\newcommand*{\eqdef}{\stackrel{\mbox{\normalfont\tiny{def}}}{=}}        
\newcommand*{\abs}[1]{\lvert #1\rvert}                
\newcommand*{\abss}[1]{\bigl\lvert #1\bigr\rvert}     
\newcommand*{\absss}[1]{\biggl\lvert #1\biggr\rvert}  

\renewcommand*{\epsilon}{\varepsilon}       
\newcommand*{\R}{\mathbb{R}}                
\newcommand*{\N}{\mathbb{N}}                
\newcommand*{\wbar}[1]{\overline{#1}}       
\newcommand*{\mcal}[1]{\mathcal{#1}}        
\newcommand*{\mbf}[1]{\mathbf{#1}}          

\DeclareMathOperator{\dist}{dist}   
\DeclareMathOperator{\ed}{ed}       
\DeclareMathOperator{\Forb}{Forb}   
\gdef\symdiff{\DOTSB\mathbin{\triangle}}    
\gdef\wj{\DOTSB\vee}        
\gdef\gj{\DOTSB\oplus}      
\gdef\crgarr{
    \mathbin{\hbox{%
    \begin{tikzpicture}[scale=0.6]%
        \draw[->] (0,0) -- (11/20,0);%
        \draw (0,1/3) -- (1/3,1/3);%
        \draw (0,0)--(0,1/3);%
    \end{tikzpicture}%
    }}
}

\begin{document}
\maketitle

\begin{abstract}
    Given a hereditary property $\mcal H$ of graphs and some $p\in[0,1]$, the edit distance function $\ed_{\mcal H}(p)$ is (asymptotically) the maximum proportion of ``edits'' (edge-additions plus edge-deletions) necessary to transform any graph of density $p$ into a member of $\mcal H$.
    For any fixed $p\in[0,1]$, $\ed_{\mcal H}(p)$ can be computed from an object known as a colored regularity graph (CRG).
    This paper is concerned with those points $p\in[0,1]$ for which infinitely many CRGs are required to compute $\ed_{\mcal H}$ on any open interval containing $p$; such a $p$ is called an accumulation point.
    We show that, as expected, $p=0$ and $p=1$ are indeed accumulation points for some hereditary properties; we additionally determine the slope of $\ed_{\mcal H}$ at these two extreme points.
    Unexpectedly, we construct a hereditary property with an accumulation point at $p=1/4$.
    Finally, we derive a significant structural property about those CRGs which occur at accumulation points.
\end{abstract}

\section{Introduction}

A \emph{hereditary property} of graphs is a family of graphs which is closed under isomorphism and under vertex-deletion.
A \emph{principal hereditary property} is of the form $\Forb(H)$ where $H$ is some fixed graph and $\Forb(H)$ is the family of all graphs which do not contain $H$ as an induced subgraph.
In fact, all hereditary properties have the form $\bigcap_{F\in\mcal F}\Forb(F)$ where $\mcal F$ is some family of graphs.
In other words, every hereditary property is defined by some collection of forbidden induced subgraphs.
It will be convenient to compress notation and write
\[
    \Forb(\mcal F)\eqdef\bigcap_{F\in\mcal F}\Forb(F).
\]

A hereditary property $\mcal H$ is said to be \emph{non-trivial} if it contains graphs of unbounded order.
Equivalently, thanks to Ramsey's theorem, $\mcal H$ is non-trivial if $\mcal H=\Forb(\mcal F)$ for some family $\mcal F$ which does not contain both a clique and an anti-clique.
Observe that if $\mcal H$ is non-trivial, then, in fact, $\mcal H$ contains graphs of every order.
\medskip

\paragraph{The edit distance function.}
For two graphs $G,H$ on vertex set $V$, we define the \emph{normalized edit distance} between $G$ and $H$ to be
\[
    \dist(G,H)\eqdef{\abs{E(G)\symdiff E(H)}\over{\abs V\choose 2}},
\]
where $\triangle$ denotes the symmetric difference.
Informally, $\dist(G,H)$ measures the number of ``edits'' (that is, edge-additions plus edge-deletions) required to transform $G$ into $H$.

If $\mcal H$ is a non-trivial hereditary property, then we may extend this definition to discuss the edit distance between $G$ and $\mcal H$ by
\[
    \dist(G,\mcal H)\eqdef\min\bigl\{\dist(G,H):H\in\mcal H\text{ s.t.\ } V(H)=V(G)\bigr\}.
\]
Note that $\dist(G,\mcal H)$ is well-defined since $\mcal H$ is closed under isomorphism and contains graphs of every order.
Finally, the \emph{edit distance function} of $\mcal H$ is defined by
\[
    \ed_{\mcal H}(p)\eqdef\limsup_{n\to\infty}\ \max\biggl\{\dist(G,\mcal H):\abs{V(G)}=n,\ \abs{E(G)}=\biggl\lfloor p {n\choose 2}\biggr\rfloor\biggr\},\qquad\text{for }p\in[0,1].
\]
Informally, $\ed_{\mcal H}(p)$ records the maximum number of edits required to transform any (large) graph of density $p$ into a member of $\mcal H$.

The edit distance function was formally defined by Balogh--Martin~\cite{BaloghMartinEditComputation}, where it was shown that $\ed_{\mcal H}$ is concave down. In a different context, Marchant--Thomason~\cite{MarchantThomasonExtremalColors} showed that $\ed_{\mcal H}$ is continuous.
For a more detailed discussion of the edit distance function and an in-depth review of known results, see~\cite{MartinSurvey}.
\medskip

The first result in this paper concerns the behavior of the edit distance function near $p=0$ and near $p=1$.
\medskip

For a family of graphs $\mcal F$, define $\chi(\mcal F)\eqdef\min\bigl\{\chi(F) : F\in\mcal F\bigr\}$ and $\wbar\chi(\mcal F)\eqdef\min\bigl\{\wbar\chi(F) : F\in\mcal F\bigr\}$ where $\wbar\chi$ denotes the clique-cover number.
Suppose that $\mcal F$ does not contain an anti-clique (i.e.\ $\chi(\mcal F)\geq 2$).
By selecting a uniformly random $(\chi(\mcal F)-1)$-coloring of the vertices of a graph and deleting any edges induced by each color class, we observe that the resulting graph is $\mcal F$-free; thus we may bound
\[
    \ed_{\Forb(\mcal F)}(p)\leq{p\over\chi(\mcal F) - 1}.
\]
While this bound is not tight for many families $\mcal F$, any known improvements are minuscule when $p$ is sufficiently small.
We show that this is no accident:

\begin{theorem}\label{thm:smallp}
    Let $\mcal F$ be a family of graphs.
    If $\mcal F$ does not contain an anti-clique, then
    \[
        \lim_{p\to 0^+}{1\over p}\ed_{\Forb(\mcal F)}(p)={1\over\chi(\mcal F) -1}.
    \]
    Similarly, if $\mcal F$ does not contain a clique, then
    \[
        \lim_{p\to 1^-}{1\over 1-p}\ed_{\Forb(\mcal F)}(p)={1\over\wbar\chi(\mcal F) -1}.
    \]
\end{theorem}
In other words, this result determines the slope of the edit distance function at $0$ and at $1$, provided the forbidden family avoids cliques and anti-cliques.

\begin{remark}
If $\mcal F$ contains an anti-clique (resp.\ clique), then the slope of $\ed_{\Forb(\mcal F)}$ at $0$ (resp.\ $1$) is a less interesting question.
This is because, if $\mcal F$ contains an anti-clique (resp.\ clique), then \Cref{prop:basicfacts} gives that one can compute $\ed_{\Forb(\mcal F)}$ over the interval $[0,1/2]$ (resp.\ $[1/2,1]$) via a finite number of structures called $p$-core CRGs, which we describe below.
\end{remark}

\paragraph{Computing the edit distance function.}
A \emph{colored regularity graph} (CRG) $K$ is a clique, together with a partition of its vertices into black and white $V(K)=VB(K)\sqcup VW(K)$ and a partition of its edges into black, white and gray $E(K)=EB(K)\sqcup EW(K)\sqcup EG(K)$.

The term CRG was, to our knowledge, coined by Alon and Stav~\cite{AlonStavFurthestGraph} but the idea traces back throughout studies of hereditary properties using Szemer\'edi's regularity lemma~\cite{SzemRegLem} such as that by Bollob\'as and Thomason~\cite{BollobasThomason1997,BollobasThomason2000}.

For a CRG $K$ and a graph $F$, we say that $F\mapsto K$ if there is a function $\phi\colon V(F)\to V(K)$ such that
\begin{itemize}
    \item For any $uv\in E(F)$, either $\phi(u)=\phi(v)\in VB(K)$ or $\phi(u)\phi(v)\in EB(K)\cup EG(K)$.
    \item For any $uv\notin E(F)$, either $\phi(u)=\phi(v)\in VW(K)$ or $\phi(u)\phi(v)\in EW(K)\cup EG(K)$.
\end{itemize}

This notion extends naturally to families of graphs, where we say that $\mcal F\mapsto K$ if there is some $F\in\mcal F$ with $F\mapsto K$.
Otherwise, we say that $\mcal F\not\mapsto K$.
For a nontrivial hereditary property $\mcal H=\Forb(\mcal F)$, denote by $\mcal K(\mcal H)$ the set of all CRGs $K$ for which $\mcal F\not\mapsto K$.
\medskip

For $p\in[0,1]$, we associate to a CRG $K$ a matrix $M_K(p)\in\R^{V(K)\times V(K)}$ defined by
\[
    \bigl(M_K(p)\bigr)_{xy}\eqdef\begin{cases}
        p & \text{if }x=y\in VW(K)\text{ or if }xy\in EW(K),\\
        1-p & \text{if }x=y\in VB(K)\text{ or if }xy\in EB(K),\\
        0 & \text{if }xy\in EG(K).
    \end{cases}
\]
From $M_K(p)$, we derive the function
\[
    g_K(p)\eqdef\min\bigl\{\langle \mu, M_K(p)\mu\rangle: \mu\in\Delta^K\bigr\},
\]
where $\Delta^K\subseteq[0,1]^{V(K)}$ is the set of all probability masses on $V(K)$ and $\langle\cdot,\cdot\rangle$ is the standard scalar product.
Note that $\Delta^K$ is naturally identified with the $(\abs{V(K)}-1)$-dimensional simplex.
\medskip

For any non-trivial hereditary property $\mcal H$ and any $p\in[0,1]$,
\begin{equation}\label{eqn:achieved}
    \ed_{\mcal H}(p)=\inf_{K\in\mcal K(\mcal H)}g_K(p)=\min_{K\in\mcal K(\mcal H)}g_K(p).
\end{equation}
Balogh--Martin~\cite{BaloghMartinEditComputation} proved the first equality and Marchant--Thomason~\cite{MarchantThomasonExtremalColors} proved the second.
In particular, for any fixed $p\in[0,1]$, the edit distance function of $\mcal H$ is determined by a single CRG.

\paragraph{Accumulation points of the edit distance function.}
Inspired by the known edit distance functions, Martin~\cite{MartinSurvey} conjectured the following extension of \cref{eqn:achieved}:
\begin{conj}[Martin~{\cite[Conjecture 1]{MartinSurvey}}]\label{conj:noacc}
    For any $\epsilon>0$ and any non-trivial hereditary property $\mcal H$, there is some finite subset $\mcal K'\subseteq\mcal K(\mcal H)$ such that
    \[
        \ed_{\mcal H}(p)=\min_{K\in\mcal K'}g_K(p),\qquad\text{for all }p\in(\epsilon,1-\epsilon).
    \]
\end{conj}

In order to approach questions of this form, we introduce the following definition.

\begin{defn}[Regular points and accumulation points of the edit distance function]
    Let $\mcal H$ be a non-trivial hereditary property.
    A point $p_0\in[0,1]$ is said to be a \emph{regular point of $\mcal H$} if there is some $\epsilon>0$ and a finite subset $\mcal K'\subseteq\mcal K(\mcal H)$ such that
    \[
        \ed_{\mcal H}(p)=\min_{K\in\mcal K'}g_K(p),\qquad\text{for every }p\in(p_0-\epsilon,p_0+\epsilon)\cap[0,1].
    \]

    Conversely, $p_0$ is said to be an \emph{accumulation point of $\mcal H$} if $p_0$ is not a regular point.
    Informally, $p_0$ is an accumulation point of $\mcal H$ if infinitely many CRGs are required to determine $\ed_{\mcal H}$ in any (relatively) open interval containing $p_0$.
\end{defn}

Appealing to compactness, \Cref{conj:noacc} is equivalent to the statement: If $\mcal H$ is a non-trivial hereditary property, then every point in $(0,1)$ is a regular point of $\mcal H$.

Martin and Riasanovsky~\cite[Theorem 39]{martinrandom} proved that $\mcal H$ does not have any accumulation points in the interval $\bigl(1-\varphi^{-1},\varphi^{-1}\bigr) = (2-\varphi,\varphi-1) = (0.382\ldots,0.618\ldots)$ where $\varphi=1.618\ldots$ is the golden ratio.
However, in this paper, we show \Cref{conj:noacc} to be false:
\begin{theorem}\label{thm:otheracc}
    $p=1/4$ is an accumulation point of $\Forb\bigl(\{K_{1,4},C_5,C_6,C_7,\ldots\}\bigr)$.
\end{theorem}

Martin~\cite[Question 5]{MartinSurvey} additionally asked if $p=0$ is indeed a accumulation point of $\Forb(K_{3,3})$ and of $\Forb(K_{2,t})$ for $t\geq 9$.
We answer this question affirmatively.

\begin{theorem}\label{thm:k33}
    $p=0$ is an accumulation point of $\Forb(K_{3,3})$ and of $\Forb(K_{2,t})$ for all $t\geq 9$.
\end{theorem}

This result follows from a general classification for when $p=0$ is an accumulation point of a hereditary property (see \Cref{thm:0acc}).

Beyond this, we show that $\Forb(K_{3,3})$ has no other accumulation points.
\begin{theorem}\label{thm:ktt}
    For any $t\in\N$, $\Forb(K_{t,t})$ has no accumulation points in the interval $(0,1]$.
\end{theorem}

In order to prove the above theorem, we establish a significant structural property about those CRGs which occur at accumulation points for a general hereditary property (see \Cref{thm:crgseq}).
We hope that this result can be beneficial to future researchers; however, it requires considerable set-up, so we are unable to state this result in the introduction.

Unfortunately, we were unable to rule out the possibility that $\Forb(K_{2,t})$ has additional accumulation points, though we believe that it does not.

\paragraph{Outline of the paper.}
We begin the paper by recalling background information about the edit distance function and CRGs in \Cref{sec:prelim}.
In \Cref{sec:slope} we establish \Cref{thm:smallp}; in order to do so, we prove a variant of the Erd\H{o}s--Stone theorem.
We then prove \Cref{thm:otheracc,thm:k33,thm:ktt} in \Cref{sec:acc}, where we additionally establish a significant structural property about those CRGs which occur at accumulation points of a general hereditary property.
Finally, we conclude with a few remarks and open questions in \Cref{sec:remarks}.

\section{Preliminaries}\label{sec:prelim}

In this section, we recall necessary definitions and results about CRGs and the edit distance function.

\begin{prop}[Martin~{\cite[Theorem 10v]{MartinSymmetrization}}]\label{prop:symmetry}
    If $\mcal H=\Forb(\mcal F)$ is a non-trivial hereditary property, then for any $p\in[0,1]$,
    \[
        \ed_{\mcal H}(p)=\ed_{\wbar{\mcal H}}(1-p),\qquad\text{where }\wbar{\mcal H}=\Forb\bigl(\{\wbar F:F\in\mcal F\}\bigr).
    \]
\end{prop}
This, in particular, allows us to restrict our attention to the interval $[0,1/2]$, since any results in this interval may be immediately translated to a result on $[0,1]$.
This is a common approach in the literature, and we employ it throughout the remainder of this paper.

\begin{defn}
    Let $K$ be a CRG and let $p\in(0,1)$.
    $K$ is said to be a \emph{$p$-core CRG} if
    \[
        g_{K'}(p)>g_K(p)
    \]
    for every proper sub-CRG $K'\subset K$.

    Equivalently, $K$ is a $p$-core CRG if there is a unique $\mu\in\Delta^K$ which achieves $g_K(p)$ and this $\mu$ has full support~\cite[Proposition 12]{MartinSymmetrization}.
\end{defn}

Note that any CRG $K$ has a $p$-core sub-CRG $L$ satisfying $g_L(p)=g_K(p)$.
\medskip

The above definition makes sense for $p\in\{0,1\}$ as well, yet it is not useful since the edit distance function is easy to compute for these values.
We will therefore define $p$-core CRGs for $p\in\{0,1\}$ structurally; these definitions are justified by a limiting argument (see \Cref{rem:0core}).
Explicitly, \Cref{prop:basicfacts} shows that the definition of $0$-core characterizes $p$-core CRGs for all $p\in (0,1/2]$ and that the definition of $1$-core characterizes $p$-core CRGs for all $p\in [1/2,1)$.
Note that \Cref{prop:basicfacts} also implies that all edges of a $1/2$-core CRG are gray.

\begin{defn}
    Let $K$ be a CRG.
    \begin{itemize}
        \item We say that $K$ is $0$-core if $K$ consists only of white and gray edges and $EW(K)\subseteq{VB(K)\choose 2}$.
        \item We say that $K$ is $1$-core if $K$ consists only of black and gray edges and $EB(K)\subseteq{VW(K)\choose 2}$.
    \end{itemize}
\end{defn}

\begin{prop}\label{prop:basicfacts}
    Fix $p\in(0,1)$ and let $K$ be a $p$-core CRG.
    \begin{enumerate}
        \item\label{basic:0core} {\normalfont(Marchant--Thomason~\cite[Lemma 3.23]{MarchantThomasonExtremalColors})} If $p\in(0,1/2]$, then $K$ is $0$-core.
        \item\label{basic:1core} {\normalfont(Marchant--Thomason~\cite[Lemma 3.23]{MarchantThomasonExtremalColors})} If $p\in[1/2,1)$, then $K$ is $1$-core.
        \item\label{basic:graydeg} {\normalfont(Martin~\cite[Lemma 15]{MartinSymmetrization})} Suppose that $p\in(0,1/2]$ and that $\mu\in\Delta^K$ achieves $g_K(p)$.
            For every $u\in VW(K)$, we have
            \[
                \mu(u)={g_K(p)\over p}.
            \]
            Furthermore, for every $u\in VB(K)$, the weighted-gray-degree of $u$ satisfies
            \[
                d_G(u)\eqdef\sum_{v\in V(K):\ uv\in EG(K)}\mu(v)={p-g_K(p)\over p}+{1-2p\over p}\mu(u).
            \]
    \end{enumerate}
\end{prop}
We remark that $d_G(u)\leq 1-\mu(u)$ for any $u\in V(K)$ since $uu\not\in EG(K)$ for any vertex $u$.
Additionally, there is an analogue to \cref{basic:graydeg} in the case that $p\in[1/2,1)$, though we will not require it thanks to \Cref{prop:symmetry}.

\begin{remark}\label{rem:0core}
    If $K$ is a $0$-core CRG, then it can be shown that there is some $p_0>0$ such that $K$ is $p$-core for all $0<p<p_0$.
    Similarly, if $K$ is a $1$-core CRG, then there is some $p_1>0$ such that $K$ is $p$-core for all $1-p_1<p<1$.

    This observation is not used to prove any of the results in this paper, so we leave the proof to the reader.
    We simply include this remark to further justify the definition of $0$-core and $1$-core CRGs.
\end{remark}

Recall that for a non-trivial hereditary property $\mcal H=\Forb(\mcal F)$, we used $\mcal K(\mcal H)$ to denote the set of all CRGs $K$ for which $\mcal F\not\mapsto K$.
We next introduce analogous notation for the CRGs in $\mcal K(\mcal H)$ which are also $p$-core for a given $p$.
\begin{defn}
    Let $\mcal H=\Forb(\mcal F)$ be a non-trivial hereditary property.
    For $p\in[0,1]$, we denote by $\mcal K_p(\mcal H)$ the set of all $p$-core CRGs $K$ for which $\mcal F\not\mapsto K$.
\end{defn}

The most important class of CRGs are those consisting only of gray edges.

\begin{defn}
    For non-negative integers $w,b$, the CRG $K(w,b)$ consists of $w$ white vertices and $b$ black vertices and has all gray edges.
\end{defn}

\begin{prop}[Martin~{\cite[Theorem 7]{MartinSymmetrization}}]\label{prop:join}
    Let $K$ and $L$ be CRGs. If $K\gj L$ denotes the CRG formed by connecting $K$ and $L$ with all gray edges, then
    \[
        {1\over g_{K\gj L}(p)}={1\over g_K(p)}+{1\over g_L(p)}.
    \]
    In particular:
    \begin{enumerate}
        \item\label{join:wb} For $0<p<1$,
            \[
                g_{K(w,b)}(p)=\biggl({w\over p}+{b\over 1-p}\biggr)^{-1}.
            \]
        \item\label{join:wk} Let $K$ be any $0$-core CRG.
            If $K$ has $w$ white vertices and $K'$ is the sub-CRG induced by the black vertices of $K$, then
            \[
                g_K(p)=\biggl({w\over p}+{1\over g_{K'}(p)}\biggr)^{-1},\qquad\text{for all }p\in(0,1).
            \]
    \end{enumerate}
\end{prop}

Intuitively, replacing a non-gray edge of a CRG $K$ with a gray edge can only decrease $g_K$.
This is formally summarized in the following statement.
\begin{claim}\label{claim:grayreplace}
    Fix $p\in(0,1)$, let $K$ be a $p$-core CRG and suppose that there is some $x\neq y\in V(K)$ with $xy\notin EG(K)$.
    If $K'$ is the CRG formed by re-coloring the edge $xy$ gray, then $g_K(p)>g_{K'}(p)$.
\end{claim}
\begin{proof}
    Let $\mu\in\Delta^K$ be the probability mass achieving $g_K(p)$; since $K$ is $p$-core we know that $\mu$ has full support.
    Now, it is easy to see that
    \[
        g_{K'}(p)\leq\langle \mu,M_{K'}(p)\mu\rangle=\langle \mu,M_K(p)\mu\rangle-2\bigl(M_K(p)\bigr)_{xy}\mu(x)\mu(y)<\langle\mu,M_K(p)\mu\rangle=g_K(p).\qedhere
    \]
\end{proof}

We will explicitly require the following consequence that holds for all CRGs, not just those that are $p$-core.
\begin{prop}\label{prop:allgray}
    Let $K$ be any CRG and fix $p\in(0,1)$.
    If $\abs{VW(K)}\leq w$ and $\abs{VB(K)}\leq b$,
    \[
        g_K(p)\geq g_{K(w,b)}(p).
    \]
    Furthermore, the above inequality is strict if one of the following holds:
    \begin{itemize}
        \item $\abs{VW(K)}<w$, or
        \item $\abs{VB(K)}<b$, or
        \item $K$ is $p$-core and has at least one non-gray edge.
    \end{itemize}
\end{prop}
\begin{proof}
    Comparing with \cref{join:wb} of \Cref{prop:join}, we observe that $g_{K(w,b)}(p)>g_{K(w',b')}(p)$ whenever either $w<w'$ and $b\leq b'$ or $w\leq w'$ and $b<b'$.
    Thus the full claim follows by applying \Cref{claim:grayreplace} to any $p$-core sub-CRG $L\subseteq K$ with $g_L(p)=g_K(p)$.
\end{proof}

\section{Slope of the edit distance function at zero}\label{sec:slope}

This section is dedicated to proving \Cref{thm:smallp}.
Appealing to the symmetry established in \Cref{prop:symmetry}, it suffices to show only that if $\mcal F$ is a family of graphs not containing an anti-clique, then
\[
    \lim_{p\to 0^+}{1\over p}\ed_{\Forb(\mcal F)}(p)={1\over\chi(\mcal F)-1}.
\]
Certainly $\mcal F\not\mapsto K(\chi(\mcal F)-1,0)$, so
\[
    \ed_{\Forb(\mcal F)}(p)\leq g_{K(\chi(\mcal F)-1,0)}(p)={p\over\chi(\mcal F)-1},\qquad\text{for all }p\in[0,1];
\]
thus the upper bound is clear.
Furthermore, observe that $\ed_{\Forb(\mcal F)}(p)\geq\ed_{\Forb(F)}(p)$ for any $F\in\mcal F$.
Since we can always find an $F\in\mcal F$ for which $\chi(F)=\chi(\mcal F)$, in order to prove \Cref{thm:smallp}, it suffices to show the following:
\begin{theorem}\label{thm:0lower}
    Let $H$ be a graph with $\chi=\chi(H)\geq 2$.
    For any $\epsilon>0$, there is some $p_0>0$ such that
    \[
        \ed_{\Forb(H)}(p)\geq (1-\epsilon){p\over\chi-1},\qquad\text{whenever }0<p<p_0.
    \]
\end{theorem}

The main ingredient in the proof of \Cref{thm:0lower} is the following variant of the Erd\H{o}s--Stone theorem~\cite{erdos-stone}.
Note that this result is purely graph-theoretic and does not involve the edit distance function and CRGs.

\begin{lemma}\label{lem:coloredES}
    Let $H$ be a graph, set $\chi=\chi(H)$ and fix $\epsilon>0$.
    There exists an integer $n_0=n_0(H,\epsilon)$ and a number $\alpha=\alpha(H,\epsilon)>0$ such that the following holds:

    If $G$ is a graph on $n> n_0$ vertices with
    \[
        \delta(G)>\biggl(1-{1\over\chi-1}+\epsilon\biggr)n,
    \]
    and $c\colon V(G)\to X$ is any coloring with $\abs{c^{-1}(x)}<\alpha n$ for each $x\in X$, then $G$ contains a (not-necessarily-induced) copy of $H$, all of whose vertices are colored differently under $c$.
\end{lemma}

\begin{proof}
    Let $K_r(t)$ denote the complete $r$-partite graph wherein each part has $t$ vertices.
    Since $H$ is a subgraph of $K_\chi(\abs{V(H)})$ it suffices to prove the claim only for the graphs $K_r(t)$.

    We prove by induction on $r$ that there are appropriate constants $n_0(r,t)\eqdef n_0\bigl(K_r(t),\epsilon\bigr)$ and $\alpha(r,t)\eqdef \alpha\bigl(K_r(t),\epsilon\bigr)$.
    Employing the convention that ${1\over 0}=+\infty$, we may consider the case of $r=1$ as our base case.
    Here, $K_1(t)$ is simply an anti-clique of size $t$, and so the result follows from setting $\alpha(1,t)=1/t$ and $n_0(1,t)=t$.
    \medskip

    For the induction step, fix $r\geq 2$ and set
    \begin{equation}\label{eqn:recdefs}
        T\eqdef{2t\over\epsilon},\quad n_0(r,t)\eqdef n_0(r-1,T),\quad\text{and}\quad \alpha(r,t)\eqdef\min\biggl\{\alpha(r-1,T),\ {\epsilon\over 2rt}{T\choose t}^{-(r-1)t}\biggr\}.
    \end{equation}
    Note that by making $\epsilon$ slightly smaller if necessary, we may suppose that $T$ is an integer; hence the above quantities are well-defined.

    Let $G=(V,E)$ be a graph on $n>n_0(r,t)$ vertices with $\delta(G)>\bigl(1-{1\over r-1}+\epsilon\bigr)n$, and let $c\colon V\to X$ be any coloring with $\abs{c^{-1}(x)}<\alpha n$ for all $x\in X$, where $\alpha=\alpha(r,t)$.
    We must show that $G$ contains a copy of $K_r(t)$, all of whose vertices are colored differently under $c$.
    From \cref{eqn:recdefs} and the induction hypothesis, we know that $G$ contains a copy of $K_{r-1}(T)$, all of whose vertices receive different colors under $c$.
    Label the parts of this $K_{r-1}(T)$ as $U_1,\dots,U_{r-1}$.

    Define the set
    \[
        A\eqdef\biggl\{(v,W_1,\dots,W_{r-1})\in V\times{U_1\choose t}\times\dots\times{U_{r-1}\choose t}:\ N(v)\supseteq\bigcup_{i=1}^{r-1}W_i\biggr\},
    \]
    where $N(v)$ denotes the neighborhood of the vertex $v$.
    We claim that we have succeeded in finding our desired copy of $K_r(t)$ if $\abs A\geq {T\choose t}^{r-1}\cdot rt\cdot\alpha n$.
    Indeed, if this is the case, then there is some $(W_1,\dots,W_{r-1})$ such that $B\eqdef\bigl\{v\in V:(v,W_1,\dots,W_{r-1})\in A\bigr\}$ satisfies $\abs B\geq rt\cdot\alpha n$.
    Since each color class of $c$ has size at most $\alpha n$, this implies that there is some $B'\subseteq B$ with $\abs{B'}\geq rt$, all of whose vertices have distinct colors under $c$.
    Removing vertices of $B'$ that share a color with a vertex in $W_1\cup\dots\cup W_{r-1}$ leaves at least $t$ vertices; these $t$ vertices along with $W_1,\dots,W_{r-1}$ forms our desired copy of $K_r(t)$.

    Therefore, the remainder of the proof is dedicated to proving that $\abs A\geq{T\choose t}^{r-1}\cdot rt\cdot\alpha n$.
    \medskip

    For $v\in V$, set $f(v)\eqdef\abss{N(v)\cap \bigcup_{i=1}^{r-1}U_i}$.
    We begin by bounding
    \begin{equation}\label{eqn:lowerfbound}
        \sum_{v\in V}f(v) = \sum_{u\in \bigcup_{i=1}^{r-1}U_i}\deg(u)> T\cdot(r-1)\cdot\biggl(1-{1\over r-1}+\epsilon\biggr)n.
    \end{equation}

    Now, let $S$ be the collection of all $v\in V$ for which $f(v)\geq T\cdot(r-1)\cdot\bigl(1-{1\over r-1}+{\epsilon\over 2}\bigr)$.
    Observe that
    \begin{align*}
        \sum_{v\in V}f(v) &\leq \abs S\cdot T\cdot(r-1)+\bigl(n-\abs S\bigr)\cdot T\cdot(r-1)\cdot\biggl(1-{1\over r-1}+{\epsilon\over 2}\biggr)\\
                          &= T\cdot(r-1)\cdot\biggl[\biggl(1-{1\over r-1}+{\epsilon\over 2}\biggr)n+\biggl({1\over r-1}-{\epsilon\over 2}\biggr)\abs S\biggr].
    \end{align*}
    Comparing with \cref{eqn:lowerfbound}, we thus bound
    \[
        \abs S>{r-1\over 2}\cdot\epsilon n\geq {\epsilon n\over 2}.
    \]

    Now, for each $v\in S$, it must be the case that
    \[
        \abs{N(v)\cap U_i}\geq{r-1\over 2}\cdot\epsilon T\geq{\epsilon T\over 2},\qquad\text{for each }i\in[r-1],
    \]
    since otherwise
    \[
        f(v)< T\cdot(r-2)+{r-1\over 2}\cdot\epsilon T=T\cdot(r-1)\cdot\biggl(1-{1\over r-1}+{\epsilon\over 2}\biggr);
    \]
    contradicting the definition of $S$.

    Since $T=2t/\epsilon$, we thus have $\abs{N(v)\cap U_i}\geq t$ for all $v\in S$ and all $i\in[r-1]$; therefore,
    \[
        \abs A\geq\abs S\geq{\epsilon n\over 2}.
    \]
    Combining this inequality with \cref{eqn:recdefs}, we finally bound
    \[
        {T\choose t}^{r-1}\cdot rt\cdot\alpha n \leq {\epsilon n\over 2}\leq \abs A,
    \]
    which concludes the proof.
\end{proof}

Returning to the edit distance function and CRGs, we next derive a simple consequence of \Cref{prop:basicfacts}.
\begin{claim}\label{claim:degrees}
    Fix $p,\epsilon$ such that $1/2>\epsilon>p>0$ and $r\in\N$ and let $K$ be a $p$-core CRG.
    Suppose that $g_K(p)<(1-\epsilon){p\over r}$ and let $\mu\in\Delta^K$ be the probability mass (with full support) achieving $g_K(p)$.
    Then for each $v\in VB(K)$,
    \[
        d_G(v)>1-{1\over r}+{\epsilon\over r},\qquad\text{and}\qquad\mu(v)<{p\over r}.
    \]
\end{claim}
\begin{proof}
    Set $g=g_K(p)$.
    From \cref{basic:graydeg} of \Cref{prop:basicfacts}, we know that
    \[
        d_G(v)={p-g\over p}+{1-2p\over p}\mu(v) > {p-g\over p}>1-{1-\epsilon\over r},
    \]
    which establishes the first piece of the claim.
    For the second part of the claim, since $d_G(v)\leq 1-\mu(v)$ and $p<\epsilon$, we obtain
    \[
        {p-g\over p}+{1-2p\over p}\mu(v)\leq 1-\mu(v)\quad\implies\quad \mu(v)\leq {g\over 1-p}<{1-\epsilon\over 1-p}\cdot{p\over r}<{p\over r}.\qedhere
    \]
\end{proof}

We are now ready to prove \Cref{thm:0lower}.

\begin{proof}[Proof of \Cref{thm:0lower}]
    Suppose the claim is false and let $H$ be a counterexample with $\chi=\chi(H)\geq 2$ minimum.
    Since $H$ is a counterexample, there is some $\epsilon>0$ and some monotone sequence $p_n\to 0$ such that, for each $n\in\N$, there is a $p_n$-core CRG $\Gamma_n$ for which
    \begin{proplist}
        \item\label{prop:notembed} $H\not\mapsto \Gamma_n$, and
        \item\label{prop:strictbound} $\displaystyle \ed_{\Forb(H)}(p_n)=g_{\Gamma_n}(p_n)<(1-\epsilon){p_n\over \chi-1}$.
    \end{proplist}

    By passing to a subsequence if necessary, we can ensure the following:
    \begin{claim}
        Each $\Gamma_n$ has only black vertices.
    \end{claim}
    \begin{proof}
        First, we know that $\abs{VW(\Gamma_n)}<\chi(H)$ for each $n$, so by passing to a subsequence if necessary, we may suppose that $\abs{VW(\Gamma_n)}=w$ for each $n\in\N$.
        Suppose for the sake of contradiction that $w>0$ and let $\Gamma_n'$ be the sub-CRG induced by the black vertices of $\Gamma_n$.
        Note that $\Gamma_n'$ is nonempty since otherwise $\Gamma_n=K(w,0)$ which has $g_{\Gamma_n}(p_n)={p_n\over w}\geq{p_n\over\chi-1}$; contradicting \cref{prop:strictbound}.
        Applying \cref{join:wk} of \Cref{prop:join}, we bound
        \[
            {1\over g_{\Gamma_n}(p_n)}={w\over p_n}+{1\over g_{\Gamma_n'}(p_n)}\quad\implies\quad g_{\Gamma_n'}(p_n)<(1-\epsilon){p_n\over \chi-w-1+\epsilon w}.
        \]

        Now, if $w=\chi-1$, then $g_{\Gamma_n'}(p_n)<p_n/\epsilon$.
        On the other hand, certainly $\abs{\Gamma_n'}< \abs{V(H)}$ or else $H\mapsto \Gamma_n$, so, through \Cref{prop:allgray},
        \[
            g_{\Gamma_n'}(p_n)\geq g_{K(0,\abs{V(H)})}(p_n)={1-p_n\over \abs{V(H)}}>{p_n\over\epsilon}
        \]
        for $n$ sufficiently large; a contradiction.

        We must therefore have $1\leq w\leq \chi -2$.
        Let $H'$ be any induced subgraph of $H$ which has $\chi(H')=\chi - w$; we know that $2\leq \chi(H')\leq\chi-1$.
        Now, it must be the case that $H'\not\mapsto \Gamma_n'$ for every $n$, or else $H\mapsto \Gamma_n$.
        Thus, for every $n$,
        \[
            \ed_{\Forb(H')}(p_n)\leq g_{\Gamma_n'}(p_n)<(1-\epsilon){p_n\over \chi-w-1+\epsilon w}<(1-\epsilon){p_n\over \chi(H')-1};
        \]
        contradicting the minimality of the graph $H$.
    \end{proof}

    We now know that each $\Gamma_n$ consists only of black vertices.
    Let $\mu_n$ denote the probability mass achieving $g_{\Gamma_n}(p_n)$.
    Fix $n$ and $N_0$ sufficiently large to be chosen later; we will consider $n$ to be large, but fixed, whereas we will consider $N_0\to\infty$.
    Define the graph $G$, which has vertex partition $V(G)=\bigsqcup_{x\in V(\Gamma_n)}V_x$ where $\abs{V_x}=\lfloor\mu_n(x)\cdot N_0\rfloor$, and whose edge-set consists of all edges between $V_x$ and $V_y$ whenever $xy\in EG(\Gamma_n)$.
    In other words, $G$ is formed by retaining only the gray edges of $\Gamma_n$ and then blowing up each vertex $x\in V(\Gamma_n)$ into an independent set of size $\lfloor\mu_n(x)\cdot N_0\rfloor$.
    Also, let $c\colon V(G)\to V(\Gamma_n)$ be the coloring defined by $c(v)=x$ if $v\in V_x$.
    We seek to show that $G$ has large minimum degree and that each color class of $c$ is sufficiently small so that we may apply \Cref{lem:coloredES}.
    \medskip

    We begin with a few observations about the structure of $G$.

    Firstly, setting $N\eqdef\abs{V(G)}$, we certainly have
    \[
        N=\sum_{x\in V(\Gamma_n)}\lfloor\mu_n(x)\cdot N_0\rfloor\quad\implies\quad N_0-\abs{V(\Gamma_n)}< N\leq N_0.
    \]

    Next, consider any $v\in V(G)$ and suppose that $c(v)=x$.
    By applying \Cref{claim:degrees} and \cref{prop:strictbound}, we bound
    \begin{align*}
        \deg(v) &= \sum_{y\in V(\Gamma_n):\ xy\in EG(\Gamma_n)}\lfloor\mu_n(x)\cdot N_0\rfloor\geq d_G(x)\cdot N-\abs{V(\Gamma_n)}\\
                &\geq\biggl(1-{1\over\chi-1}+{\epsilon\over\chi-1}\biggr)N-\abs{V(\Gamma_n)}.
    \end{align*}

    Finally, by selecting $n$ sufficiently large so that $p_n<\epsilon$, for any $x\in V(\Gamma_n)$
    \[
        \abs{c^{-1}(x)}=\lfloor\mu_n(x)\cdot N_0\rfloor\leq \mu_n(v)\cdot N+\abs{V(\Gamma_n)}<{p_n\over\chi-1}\cdot N+\abs{V(\Gamma_n)},
    \]
    thanks again to \Cref{claim:degrees} and \cref{prop:strictbound}.
    \medskip

    Thus, by first selecting $n$ sufficiently large (so that $p_n$ is sufficiently small) and then selecting $N_0$ (and hence $N$) sufficiently large compared to $\abs{V(\Gamma_n)}$, we may apply \Cref{lem:coloredES} to the graph $G$ to find a copy of $H$, all of whose vertices are colored differently under $c$.

    Suppose that $\phi\colon V(H)\to V(G)$ realizes $H$ as such a subgraph of $G$ and consider the map $c\circ\phi\colon V(H)\to V(\Gamma_n)$.
    Since the vertices of this copy of $H$ are all colored differently under $c$, we know that $c\circ\phi$ is an injection.
    We therefore find that $c\circ\phi$ maps edges of $H$ to gray edges of $\Gamma_n$ and maps non-edges of $H$ to either white or gray edges of $\Gamma_n$.
    This, however, means that $H\mapsto \Gamma_n$; contradicting \cref{prop:notembed}.
\end{proof}

\section{Accumulation points}\label{sec:acc}

We begin by establishing a necessary and sufficient condition for $0$ to be an accumulation point of a hereditary property.
Thanks to \Cref{prop:symmetry}, there is an analogous classification for when $1$ is an accumulation point.

\begin{theorem}\label{thm:0acc}
    Let $\mcal H=\Forb(\mcal F)$ be a non-trivial hereditary property and set $\chi=\chi(\mcal F)$.
    For $p\in(0,1/2)$, define
    \[
        q_{\mcal H}(p)\eqdef\min\bigl\{ g_K(p): \abs{VW(K)}=\chi-1,\ K\in\mcal K_0(\mcal H)\bigr\}.
    \]
    $p=0$ is a regular point of $\mcal H$ if and only if
    \[
        \ed_{\mcal H}(p)=q_{\mcal H}(p),\qquad\text{for all $p$ sufficiently small.}
    \]
\end{theorem}
\begin{proof}
    Recalling that $\mcal K_0(\mcal H)$ is the set of $0$-core CRGs in $\mcal K(\mcal H)$, we observe that $\bigl\{K\in\mcal K_0(\mcal H):\abs{VW(K)}=\chi-1\bigr\}$ is finite.
    Indeed, if $F\in\mcal F$ has $\chi(F)=\chi$, then $F\mapsto K$ whenever $K$ is a $0$-core CRG, $\abs{VW(K)}\geq\chi-1$ and $\abs{VB(K)}\geq\abs{V(F)}$.
    Therefore $q_{\mcal H}(p)$ is well-defined and also the ``if'' direction is immediate.
    \medskip

    To establish the ``only if'' direction, suppose that $p=0$ is a regular point of $\mcal H$; thus, we can find some finite family $\mcal K'\subseteq\mcal K(\mcal H)$ with the property that
    \[
        \ed_{\mcal H}(p)=\min_{K\in\mcal K'}g_K(p)
    \]
    for all $p$ sufficiently small.
    Certainly we may assume that $\mcal K'\subseteq\mcal K_0(\mcal H)$.

    Set
    \[
        N\eqdef\max_{K\in\mcal K'}\abs{VB(K)},
    \]
    and consider some $K\in\mcal K'$.
    If $\abs{VW(K)}=\chi-1$, then $g_K(p)\geq q_{\mcal H}(p)$ by definition, so suppose that $\abs{VW(K)}=w\leq\chi-2$.
    In this case, we apply \Cref{prop:allgray} to bound
    \[
        g_K(p)\geq g_{K(w,N)}(p)={(1-p)p\over (1-p)w+pN}>{p\over\chi -1}=g_{K(\chi-1,0)}(p).
    \]
    provided $p<1/(N+1)$.
    Since certainly $K(\chi-1,0)\in\mcal K_0(\mcal H)$, we have shown that, for all $p$ sufficiently small,
    \[
        q_{\mcal H}(p)\geq\ed_{\mcal H}(p)=\min_{K\in\mcal K'}g_K(p)\geq q_{\mcal H}(p).\qedhere
    \]
\end{proof}

Using \Cref{thm:0acc}, we can quickly verify that $p=0$ is indeed an accumulation point for some hereditary properties.

\begin{proof}[Proof of \Cref{thm:k33}]
    First, set $\mcal H=\Forb(K_{3,3})$.
    If $K\in\mcal K_0(\mcal H)$ has exactly one white vertex, then $K$ has at most two black vertices.
    Since $K_{3,3}\not\mapsto K(1,2)$, an appeal to \Cref{prop:allgray} implies that
    \[
        q_{\mcal H}(p)={p(1-p)\over 1+p},\qquad\text{for every }p\in(0,1/2).
    \]
    However, Marchant and Thomason~\cite[Example 5.19]{MarchantThomasonExtremalColors} showed that there is a monotone sequence $p_n\to 0$ such that
    \[
        \ed_{\mcal H}(p_n)<{p_n(1-p_n)\over 1+p_n},\qquad\text{for every }n\in\N.
    \]
    Thus, thanks to \Cref{thm:0acc}, we know that $p=0$ is indeed an accumulation point of $\mcal H$.
    \medskip

    Next set $\mcal H=\Forb(K_{2,t})$ for any $t\geq 9$.
    Observe that if $K\in\mcal K_0(\mcal H)$ has exactly one white vertex, then $K$ has at most one black vertex.
    Since $K_{2,t}\not\mapsto K(1,1)$ for $t\geq 2$, another appeal to \Cref{prop:allgray} implies that
    \[
        q_{\mcal H}(p)=p(1-p),\qquad\text{for every }p\in(0,1/2).
    \]
    However, Martin and McKay~\cite[Theorem 8 \& Corollary 9]{martink2t} showed that there is a monotone sequence $p_n\to 0$ such that
    \[
        \ed_{\mcal H}(p_n)<p_n(1-p_n),\qquad\text{for every }n\in\N.
    \]
    Thus, again thanks to \Cref{thm:0acc}, we know that $p=0$ is indeed an accumulation point of $\mcal H$.
\end{proof}

We prove that $p=0$ is the \emph{only} accumulation point of $\Forb(K_{3,3})$ in the proof of \Cref{thm:ktt} in \Cref{sec:k33}.
\medskip

Before we get to this, we prove \Cref{thm:otheracc}.
To do so, we need an intermediate result.
Firstly, for two CRGs $K,L$, denote by $K\wj L$ the CRG formed by connecting $K$ and $L$ with all-white edges.
Next, denote by $\mcal P_n$ the $n$-vertex CRG with only black vertices, whose gray edges form an $n$-vertex path and all of whose other edges are white.

\begin{lemma}\label{lem:pathbound}
    For any $n_1,\dots,n_\ell\in\N$, there is some $\epsilon>0$ such that
    \[
        g_{\mcal P_{n_1}\wj\dots\wj\mcal P_{n_\ell}}(p)>p,\qquad\text{for all }p\in(1/4-\epsilon,1/4+\epsilon)
    \]
\end{lemma}
\begin{proof}
    Since $g_K$ is a continuous function for any CRG $K$, it suffices to show only that
    \[
        g_{\mcal P_{n_1}\wj\dots\wj\mcal P_{n_\ell}}(1/4)>1/4.
    \]
    Since any sub-CRG of $\mcal P_{n_1}\wj\dots\wj\mcal P_{n_\ell}$ also has the form $\mcal P_{m_1}\wj\dots\wj \mcal P_{m_t}$ for some $m_1,\dots,m_t\in\N$, we may assume, without loss of generality, that $K\eqdef\mcal P_{n_1}\wj\dots\wj\mcal P_{n_\ell}$ is $1/4$-core.

    For a vertex $v\in V(K)$, denote by $NG(v)\eqdef\{u\in V(K):uv\in EG(K)\}$ the gray-neighborhood of $v$.
    Set $g=g_K(1/4)$ and let $\mu\in\Delta^K$ be the probability mass achieving $g$.

    By applying \cref{basic:graydeg} of \Cref{prop:basicfacts}, we know that
    \[
        2\mu(v)+1-4g=\sum_{u\in NG(v)}\mu(u),\qquad\text{for all }v\in V(K).
    \]
    Supposing that $\abs{V(K)}=n$, summing both sides over all $v\in V(K)$ then yields
    \begin{align*}
        2+n(1-4g) &= \sum_{v\in V(K)}\sum_{u\in NG(v)}\mu(u)=\sum_{v\in V(K)}\abs{NG(v)}\cdot\mu(v)\\
                  &= 2-\sum_{v:\abs{NG(v)}=1}\mu(v)-2\sum_{v:\abs{NG(v)}=0}\mu(v)\\
        \implies g &= {1\over 4}+{1\over 4n}\sum_{v:\abs{NG(v)}=1}\mu(v)+{1\over 2n}\sum_{v:\abs{NG(v)}=0}\mu(v)>{1\over 4},
    \end{align*}
    since $\mu$ has full support and $K$ certainly has at least one vertex $v\in V(K)$ with $\abs{NG(v)}\in\{0,1\}$.
\end{proof}

We now demonstrate a hereditary property which has an accumulation point at $p=1/4$.

\begin{proof}[Proof of \Cref{thm:otheracc}]
    Set $\mcal H=\Forb(\mcal F)$ where $\mcal F=\{K_{1,4},C_5,C_6,\ldots\}$.
    We begin by determining $\mcal K_0(\mcal H)$.
    \begin{claim}
        \[
            \mcal K_0(\mcal H)=\bigl\{K(1,0)\bigr\}\cup\bigl\{\mcal P_{n_1}\wj\dots\wj\mcal P_{n_\ell}:n_1,\dots,n_\ell\in\N\bigr\}.
        \]
    \end{claim}
    \begin{proof}
        Since $\mcal F$ does not contain an anti-clique, we know that $K(1,0)\in\mcal K_0(\mcal H)$.
        Furthermore, since $K_{1,4}\mapsto K(2,0)$ and $K_{1,4}\mapsto K(1,1)$, we observe that $K(1,0)$ is the only member of $\mcal K_0(\mcal H)$ which has at least one white vertex.

        Next it is easy to observe that $K_{1,4}\not\mapsto\mcal P_{n_1}\wj\dots\wj \mcal P_{n_\ell}$ and that if $C_m\mapsto\mcal P_{n_1}\wj\dots\wj\mcal P_{n_\ell}$, then $m\leq 4$.
        Thus, all that is left to show is that if $K\in\mcal K_0(\mcal H)$ has only black vertices, then $K=\mcal P_{n_1}\wj\dots\wj\mcal P_{n_\ell}$ for some $n_1,\dots,n_\ell\in\N$.

        Again, since $K_{1,4}\not\mapsto K$, we observe that every vertex of $K$ is incident to at most two gray edges.
        Beyond this, if $K$ contains an $n$-vertex gray-edge cycle, then $C_m\mapsto K$ for every $n\leq m\leq 2n$.
        We conclude that the gray-edges of $K$ form a linear forest and so we have established the claim.
    \end{proof}

    We are now ready to show that $p=1/4$ is an accumulation point of $\mcal H$.
    Suppose to the contrary that there exists some $\epsilon>0$ and some finite family $\mcal K'\subseteq\mcal K(\mcal H)$ such that
    \[
        \ed_{\mcal H}(p)=\min_{K\in\mcal K'}g_K(p),\qquad\text{for all }p\in[1/4,1/4+\epsilon).
    \]
    Certainly we may suppose that $\epsilon<1/4$ and that $\mcal K'\subseteq\mcal K_0(\mcal H)$.

    Now, we know that $g_{K(1,0)}(p)=p$ for all $p\in[0,1]$, so by additionally applying \Cref{lem:pathbound} and using the fact that $\mcal K'$ is finite, we observe that there is some $\epsilon_1>0$ for which
    \begin{equation}\label{eqn:lowerfor14}
        \ed_{\mcal H}(p)=\min_{K\in\mcal K'}g_K(p)\geq p,\qquad\text{for all }p\in[1/4,1/4+\epsilon_1).
    \end{equation}

    However, for $n\geq 3$, letting $\mbf u\in\Delta^{\mcal P_n}$ denote the uniform distribution, we compute
    \[
        g_{\mcal P_n}(p)\leq \langle\mbf u, M_{\mcal P_n}(1/4)\mbf u\rangle={1-p\over n}+p\cdot{(n-2)(n-3)+2(n-1)\over n^2}=p-{4p-1\over n}+{4p\over n^2}.
    \]
    Therefore, since $\mcal F\not\mapsto\mcal P_n$, we see that for any $p>1/4$ and $n$ sufficiently large, we have
    \[
        \ed_{\mcal H}(p)\leq g_{\mcal P_n}(p)<p,
    \]
    contradicting \cref{eqn:lowerfor14}.
\end{proof}

\begin{remark}
    The above proof actually establishes a slightly stronger claim; namely, if $\mcal C$ is any collection of cycles satisfying
    \begin{itemize}
        \item If $C_m\in\mcal C$, then $m\geq 5$, and
        \item For every $n\geq 3$, there is some $n\leq m\leq 2n$ for which $C_m\in\mcal C$,
    \end{itemize}
    then $\Forb(\{K_{1,4}\}\cup\mcal C)$ has an accumulation point at $p=1/4$.
\end{remark}

\subsection{CRGs at accumulation points}\label{sec:k33}
In this section, we establish structural properties of CRGs at accumulation points of the edit distance function.

In order to do so, we will require the framework established by Marchant and Thomason~\cite{MarchantThomasonExtremalColors}.
In their paper, Marchant and Thomason work with a function different from the edit distance function, which they show to be equivalent to $1-\ed_{\mcal H}$ (see~\cite[Section 2.2]{MarchantThomasonExtremalColors}).
Moreover, they use somewhat different notation. Below we translate their framework to our situation of edit distance, and we modify their notation to be more consistent with our own.
\medskip

We first define a class of objects closely related to CRGs.
\begin{defn}
    A \emph{colored graph} $G$ is a clique with a coloring of its edges using the colors black, white and gray: $E(G)=EB(G)\sqcup EW(G)\sqcup EG(G)$. We use $\abs G$ to denote the number of vertices of a colored graph $G$.

    For a colored graph $G$ and a number $p\in[0,1]$, define the following quantities:
    \begin{itemize}
        \item For an edge $e\in E(G)$, the $p$-weight of $e$ is defined as
            \[
                w_p(e)\eqdef\begin{cases}
                    p & \text{if }e\in EW(G),\\
                    1-p & \text{if }e\in EB(G),\\
                    0 & \text{if }e\in EG(G).
                \end{cases}
            \]
        \item For $v\in V(G)$, the $p$-degree of $v$ is defined as
            \[
                d_p(v)\eqdef\sum_{u\in V(G)}w_p(uv).
            \]
        \item Finally, the maximum-$p$-degree of $G$ is defined as
            \[
                \Delta_p(G)\eqdef\max_{v\in V(G)}d_p(v).
            \]
    \end{itemize}
\end{defn}

The only difference between a colored graph and a CRG is that a colored graph does not come equipped with a coloring of its vertices.

\begin{defn}
For two colored graphs $G,H$, we say that $G\sqsubseteq H$ if there is an injection $\phi\colon V(G)\to V(H)$ so that
\begin{itemize}
    \item If $uv\in EB(G)$, then $\phi(u)\phi(v)\in EB(H)\cup EG(H)$, and
    \item If $uv\in EW(G)$, then $\phi(u)\phi(v)\in EW(H)\cup EG(H)$, and
    \item If $uv\in EG(G)$, then $\phi(u)\phi(v)\in EG(H)$.
\end{itemize}
\end{defn}
In other words, $G\sqsubseteq H$ if $G$ is formed from $H$ by deleting some vertices and recoloring some gray edges with either white or black.
Note that $\sqsubseteq$ is a partial-order on the set of colored graphs and that $G\subseteq H\implies G\sqsubseteq H$.
\medskip

A single colored graph is uninteresting on its own --- we instead need to work with sequences of colored graphs.
\begin{defn}
    A \emph{colored graph sequence} is a sequence of colored graphs $(G_n)_{n\in\N}$ of unbounded order.
\end{defn}
It is useful to think of a colored graph sequence as being an approximation of a CRG, though we will see later that colored graph sequences can encode more intricate information.

From a single CRG, we can generate a colored graph sequence by ``blowing-up'' the vertices into colored cliques:
\begin{defn}\label{defn:blowup}
    Let $K$ be a CRG and let $m\in\N$.
    The colored graph $m\times K$ has vertex-set $V=\bigsqcup_{x\in V(K)}V_x$ where $\abs{V_x}=m$ for every $x\in V(K)$.
    Furthermore, all edges between vertices in $V_x$ have the same color as the vertex $x\in V(K)$, and all edges between $V_x$ and $V_y$ have the same color as the edge $xy\in E(K)$.
\end{defn}

In the opposite direction, the following definition lays out a template for passing from a colored graph sequence to a CRG.

\begin{defn}
    Let $(G_n)_{n\in\N}$ be a colored graph sequence and let $K$ be a CRG.
    We say that $K\crgarr (G_n)$ if for every $m\in\N$, there is some $n_m\in\N$ for which $m\times K\sqsubseteq G_{n_m}$.
\end{defn}
Observe that if $K\crgarr(G_n)$, then for any $m\in\N$, there are, in fact, infinitely many $n$'s for which $m\times K\sqsubseteq G_n$.
Furthermore, certainly if $K\subseteq L$ and $L\crgarr(G_n)$, then also $K\crgarr(G_n)$.

To motivate the above definition, consider a graph $H$; we may identify $H$ with a colored graph by coloring the edges of $H$ black and coloring the non-edges of $H$ white.
Observe that for a CRG $K$, if $H\mapsto K$, then $H\sqsubseteq m\times K$ for $m$ sufficiently large (namely, $m\geq\abs{V(H)}$).
Therefore, if $(G_n)$ is a colored graph sequence with the property that $H\not\sqsubseteq G_n$ for all $n$, then for any CRG $K$ with $K\crgarr (G_n)$, it must be the case that $H\not\mapsto K$.
\medskip

We will require the following results of Marchant and Thomason from~\cite{MarchantThomasonExtremalColors}.
First, a CRG $L$ is said to be an \emph{extension} of a CRG $K$ if $K$ is obtained by deleting a single vertex from $L$.

\begin{lemma}[Marchant--Thomason~{\cite[Lemma 3.11]{MarchantThomasonExtremalColors}}]\label{mt:extension}
    Let $(G_n)_{n\in\N}$ be a colored graph sequence and suppose that $K$ is a CRG for which $K\crgarr(G_n)$.
    Then for any $\mu\in\Delta^K$, there is an extension $L$ of $K$ for which $L\crgarr(G_n)$ and the vertex $\{v\}=V(L)\setminus V(K)$ satisfies
    \[
        \sum_{u\in V(K)}\mu(u)w_p(uv)\leq\limsup_{n\to\infty}{\Delta_p(G_n)\over\abs{G_n}}.
    \]
\end{lemma}

\begin{lemma}[Marchant--Thomason~{\cite[Lemmas 3.13 \& 3.22]{MarchantThomasonExtremalColors}}]\label{mt:betterext}
    Fix $p\in(0,1/2)$, let $K$ be a $p$-core CRG and let $\mu\in\Delta^K$ be the probability mass achieving $g_K(p)$.
    Suppose that $L$ is an extension of $K$ such that
    \[
        \sum_{u\in V(K)}\mu(u)w_p(uv)<g_K(p),
    \]
    where $\{v\}=V(L)\setminus V(K)$.
    If $L'$ is any sub-CRG of $L$ which is $p$-core and has $g_{L'}(p)=g_L(p)$, then $g_{L'}(p)<g_K(p)$ and $\abs{VW(L')}\geq\abs{VW(K)}$.
\end{lemma}

The bulk of the arguments in this section are built on variants of the following observation, which, informally, states that a large number of black vertices of a $0$-core CRG behave as if they were a single white vertex.
\begin{prop}\label{prop:manyblack}
    Let $K$ be a $0$-core CRG and fix $n\in\N$.
    If $\abs{VW(K)}=t$ and $\abs{VB(K)}\geq n$, then $n\times K(t+1,0)\sqsubseteq n\times K$.
\end{prop}
\begin{proof}
    Let $V_1,\dots,V_{t+1}$ be the vertex partition of $V(n\times K(t+1,0))$ defined in \Cref{defn:blowup}; label $V_i=\{v_i^1,\dots,v_i^n\}$.
    Additionally, let $W_1,\dots,W_t,B_1,\dots,B_m$ be the vertex partition of $V(n\times K)$ where $W_1,\dots,W_t$ correspond to the white vertices of $K$ and $B_1,\dots,B_m$ correspond to the black vertices of $K$.
    Label $W_i=\{w_i^1,\dots,w_i^n\}$ and $B_i=\{b_i^1,\dots,b_i^n\}$.

    Define the map $\phi\colon V(n\times K(t+1,0))\to V(n\times K)$ by
    \[
        \phi(v_i^j)=\begin{cases}
            w_i^j & \text{if }i\in[t],\\
            b_j^1 & \text{if }i=t+1,
        \end{cases}
    \]
    which is well-defined and an injection since $m\geq n$.
    Since $K$ is $0$-core, all edges induced by $\phi(V_i)$ for any $i\in[t+1]$ are either white or gray, and all edges between $\phi(V_i)$ and $\phi(V_j)$ for $i\neq j\in[t+1]$ are gray.
    Therefore $n\times K(t+1,0)\sqsubseteq n\times K$.
\end{proof}

The following result was implicitly established by Marchant and Thomason in the proof of Theorem 3.25 in~\cite{MarchantThomasonExtremalColors}.
\begin{lemma}\label{lem:maxwhite}
    Fix $p\in(0,1/2)$, let $(G_n)_{n\in\N}$ be a colored graph sequence and set
    \[
        w=\sup\bigl\{t\in\N\cup\{0\}:K(t,0)\crgarr(G_n)\bigr\}.
    \]
    If $w$ is finite, then there is a $p$-core CRG $K$ such that $K\crgarr(G_n)$, $K$ has at least $w$ white vertices and
    \[
        g_K(p)\leq\limsup_{n\to\infty}{\Delta_p(G_n)\over\abs{G_n}}.
    \]
\end{lemma}
\begin{proof}
    Set
    \[
        \Delta_p\eqdef\limsup_{n\to\infty}{\Delta_p(G_n)\over\abs{G_n}}.
    \]

    Suppose for the sake of contradiction that the claim is false; we aim to show that $K(w+1,0)\crgarr(G_n)$.
    Under this assumption, we begin by defining a sequence of $p$-core CRGs $(\Gamma_n)_{n\in\N}$ such that
    \begin{itemize}
        \item $\abs{VW(\Gamma_n)}\geq w$, and
        \item $\Gamma_n\crgarr(G_n)$, and
        \item $g_{\Gamma_n}(p)>g_{\Gamma_{n+1}}(p)$.
    \end{itemize}

    Begin with $\Gamma_1= K(w,0)$, which is $p$-core, has $w$ white vertices and $\Gamma_1\crgarr(G_n)$ by definition.
    Assuming $\Gamma_n$ has been defined, we define $\Gamma_{n+1}$ as follows:

    Let $L$ be an extension of $\Gamma_n$ guaranteed by \Cref{mt:extension}, so $L\crgarr(G_n)$ and
    \[
        \sum_{u\in V(\Gamma_n)}\mu(u)w_p(uv)\leq\Delta_p<g_{\Gamma_n}(p),
    \]
    where $\{v\}=V(L)\setminus V(\Gamma_n)$ and the latter inequality holds by assumption.
    Now, let $\Gamma_{n+1}$ be any $p$-core sub-CRG of $L$ with $g_{\Gamma_{n+1}}(p)=g_L(p)$.
    Since $\Gamma_{n+1}\subseteq L$ and $L\crgarr(G_n)$, we must have $\Gamma_{n+1}\crgarr(G_n)$.
    Additionally, thanks to \Cref{mt:betterext}, we know that $g_{\Gamma_{n+1}}(p)<g_{\Gamma_n}(p)$ and that $\abs{VW(\Gamma_{n+1})}\geq\abs{VW(\Gamma_n)}\geq w$.
    Thus, $\Gamma_{n+1}$ satisfies the claimed properties.
    \medskip

    We claim now that $K(w+1,0)\crgarr(G_n)$, which will contradict the definition of $w$ and conclude the proof.
    Firstly, if $\abs{VW(\Gamma_n)}>w$ for some $n\in\N$, then $K(w+1,0)\subseteq \Gamma_n$ and so $K(w+1,0)\crgarr(G_n)$.
    We may thus suppose that $\abs{VW(\Gamma_n)}=w$ for each $n\in\N$.

    Fix any $N\in\N$.
    Since $g_{\Gamma_1}(p)>g_{\Gamma_2}(p)>\cdots$, we know that the $\Gamma_n$'s are distinct; thus, since each $\Gamma_n$ has exactly $w$ white vertices, there must be some $t\in\N$ for which $\abs{VB(\Gamma_t)}\geq N$.
    Now, since $\Gamma_t\crgarr(G_n)$, there is some $t'\in\N$ for which $N\times \Gamma_t\sqsubseteq G_{t'}$.
    Thus, by applying \Cref{prop:manyblack}, we have $N\times K(w+1,0)\sqsubseteq N\times \Gamma_t\sqsubseteq G_{t'}$.
    Since $N\in\N$ was arbitrary, we conclude that $K(w+1,0)\crgarr(G_n)$, which concludes the proof.
\end{proof}

We next demonstrate a more refined construction of a colored graph sequence from a single CRG wherein we ``blow-up'' the vertices into colored cliques of various sizes.
\begin{defn}\label{defn:probblowup}
    Let $K$ be a CRG and fix any $\mu\in\Delta^K$.
    For a positive integer $n$, denote by $K[\mu,n]$ the colored graph with vertex-partition $V(K[\mu,n])=\bigsqcup_{x\in V(K)}V_x$ where $\abs{V_x}=\lfloor \mu(x)\cdot n\rfloor$ for every $x\in V(K)$.
    Furthermore, all edges between vertices in $V_x$ have the same color as the vertex $x\in V(K)$, and all edges between $V_x$ and $V_y$ have the same color as the edge $xy\in E(K)$.
\end{defn}

Notice that if $k=\abs{V(K)}$, then $m\times K=K[\mbf u,km]$ where $\mbf u \in\Delta^K$ is the uniform distribution.

Sequences of colored graphs constructed in this manner mirror the properties of the original CRG.

\begin{prop}\label{prop:blowup}
    Fix $p\in(0,1/2]$.
    If $K$ is a $p$-core CRG and $\mu\in\Delta^K$ is the probability mass achieving $g_K(p)$, then
    \[
        g_K(p)=\lim_{n\to\infty}{\Delta_p(K[\mu,n])\over \abs{K[\mu,n]}}.
    \]
\end{prop}
\begin{proof}
    For notational ease, set $G_n\eqdef K[\mu,n]$, and note that
    \[
        \abs{G_n}=\sum_{x\in V(K)}\lfloor\mu(x)\cdot n\rfloor\quad\implies\quad \abs{G_n}=n+O(1),
    \]
    for $n$ sufficiently large.

    Fix any $v\in V(G_n)$ and suppose that $v\in V_x$ where $x\in V(K)$.
    If $x\in VW(K)$, then
    \[
        d_p(v)=p\cdot\bigl(\lfloor\mu(x)\cdot n\rfloor-1\bigr).
    \]
    since $v$ is connected by white edges to every $u\in V_x$ and by gray edges to every other vertex of $G_n$.
    Thus, for $n$ sufficiently large, we have
    \[
        {d_p(v)\over\abs{G_n}}=p\cdot\mu(x)+O(1/n)=g_K(p)+O(1/n),
    \]
    thanks to \cref{basic:graydeg} of \Cref{prop:basicfacts}.

    On the other hand, if $x\in VB(K)$, then $v$ is connected by black edges to every $u\in V_x$ and by either gray or white edges to every other vertex of $G_n$, so
    \begin{align*}
        d_p(v) &=(1-p)\cdot\bigl(\lfloor\mu(x)\cdot n\rfloor-1\bigr)+p\cdot\sum_{y\in V(K):xy\in EW(K)}\lfloor \mu(y)\cdot n\rfloor\\
               &=(1-p)\cdot\bigl(\lfloor\mu(x)\cdot n\rfloor-1\bigr)+p\cdot\biggl(\abs{G_n}-\lfloor\mu(x)\cdot n\rfloor-\sum_{y\in V(K):xy\in EG(K)}\lfloor \mu(y)\cdot n\rfloor\biggr).
    \end{align*}
    Thus, by again appealing to \cref{basic:graydeg} of \Cref{prop:basicfacts}, for $n$ sufficiently large, we have
    \begin{align*}
        {d_p(v)\over\abs{G_n}} &= (1-p)\mu(x)+p\cdot\bigl(1-\mu(x)-d_G(x)\bigr)+O(1/n)\\
                               &= (1-p)\mu(x)+p\cdot\biggl(1-\mu(x)-{p-g_K(p)\over p}-{1-2p\over p}\mu(x)\biggr)+O(1/n)\\
                               &= g_K(p)+O(1/n).
    \end{align*}

    We conclude that
    \[
        {\Delta_p(G_n)\over\abs{G_n}}=g_K(p)+O(1/n),
    \]
    and so the claim follows.
\end{proof}

To go further, we consider modifying a CRG by replacing a white vertex by a collection of black vertices.
\begin{defn}\label{defn:dalmatian}
    A \emph{dalmatian CRG} of size $\ell$ is a CRG consisting of $\ell$ black vertices with all white edges.

    Let $K$ be a $0$-core CRG with $\abs{VW(K)}=w$, fix an integer $\ell\in\N$ and an integer $1\leq r\leq w$.
    The CRG $K^\ell(r)$ is defined by replacing $r$ of $K$'s white vertices by dalmatian CRGs of size $\ell$ and otherwise leaving $K$ unchanged.
\end{defn}

Dalmatian CRGs were introduced by Martin--Riasanovsky~\cite{martinrandom}.
\medskip

The following lemma expands on \Cref{prop:manyblack} and is the key step in the main result of this section.
Intuitively, it states that if a collection of black vertices in a $0$-core CRG behave a single white vertex (as per \Cref{prop:manyblack}), then they actually behave like a dalmatian CRG.
\begin{lemma}\label{lem:dalmatian}
    Let $K,L$ be CRGs where $K$ is $0$-core and fix $n\in\N$.
    Let $m$ be another positive integer and let $\mu\in\Delta^L$ be any probability mass satisfying $m\cdot\mu(x)\geq n$ for all $x\in VB(L)$.

    If $\abs{VW(K)}=s>t=\abs{VW(L)}$ and $n\times K\sqsubseteq L[\mu,m]$, then also $n\times K^n(r)\sqsubseteq L[\mu,m]$ for any $1\leq r\leq s-t$.
\end{lemma}
\begin{proof}
    Observe that if $r>1$, then $K^n(r)=\bigl(K^n(r-1)\bigr)^n(1)$.
    Furthermore, $K^n(1)$ is still $0$-core and $\abs{VW(K^n(1))}=\abs{VW(K)}-1$.
    Thus, it suffices to establish the claim only for $r=1$ since then the full claim follows by a straight-forward induction.
    \medskip

    Let $W_1,\dots,W_s,B_1,\dots,B_{k}$ be the vertex partition of the colored graph $n\times K$ given by \Cref{defn:blowup}, where the $W_i$'s correspond to the white vertices of $K$ and the $B_i$'s correspond to the black vertices of $K$.
    Similarly, let $W_1',\dots,W_t',B_1',\dots,B_{k'}'$ be the vertex partition of the colored graph $L[\mu,m]$ given by \Cref{defn:probblowup}.
    Furthermore, let $\phi\colon V(n\times K)\to V(L[\mu,m])$ be an injection which exhibits $n\times K\sqsubseteq L[\mu,m]$.
    \medskip

    Fix any $i\neq j\in[s]$ and $\ell\in[t]$; since each edge between $W_i$ and $W_j$ is gray and all edges inside $W_\ell'$ are white, it cannot be the case that both $\phi(W_i)\cap W_\ell'$ and $\phi(W_j)\cap W_\ell'$ are nonempty.
    Since $s>t$, we may therefore suppose, without loss of generality, that $\phi(W_s)\subseteq\bigcup_{i=1}^{k'}B_i'$.
    Similarly, all edges within $W_s$ are white and all edges within each $B_i'$ are black, so no two vertices of $W_s$ can be mapped to the same $B_i'$.
    Thus, labeling $W_s=\{w_1,\dots,w_n\}$, and relabeling the $B_i'$'s if necessary, we may suppose that $\phi(w_i)\in B_i'$.

    Now, again since $K$ is $0$-core, every edge incident to $w_i$ is white or gray.
    Furthermore, since all edges within each $B_i'$ are black, we observe that $\phi^{-1}(B_i')=\{w_i\}$ for all $i\in[n]$.
    Additionally, since $\phi$ exhibits $n\times K\sqsubseteq L[\mu,m]$, we know that all edges between $B_i'$ and $\phi(v)$ for any $i\in[n]$ and any $v\in V(n\times K)\setminus W_s$ must be gray.
    \medskip

    We now turn our attention to $n\times K^n(1)$, which we can consider to have vertex-partition
    \[
        W_1,\dots,W_{s-1},D_1,\dots,D_n,B_1,\dots,B_k,
    \]
    where $D_1,\dots,D_n$ correspond to the black vertices in the new dalmatian set.
    Consider a function $\phi'\colon V(n\times K^n(1))\to V(L[\mu,m])$ with the following properties:
    \begin{itemize}
        \item $\phi'(v)=\phi(v)$ for all $v\in W_1\cup\dots\cup W_{s-1}\cup B_1\cup\dots\cup B_k$, and
        \item $\phi'(d)\in B_i'$ for all $i\in[n]$ and all $d\in D_i$.
    \end{itemize}
    By the assumption on $m$ and $\mu$, for all $i\in[n]$, we have $\abs{B_i'}\geq n=\abs{D_i}$, so such a $\phi'$ exists which is also an injection.
    Based on the properties of $\phi$ laid out above, this $\phi'$ realizes $n\times K^n(1)\sqsubseteq L[\mu,m]$.
\end{proof}

We are now ready to state and prove the main result in this section.

\begin{theorem}\label{thm:crgseq}
    Let $\mcal H=\Forb(\mcal F)$ be a non-trivial hereditary property, fix $p\in(0,1/2)$ and suppose that there is a sequence of CRGs $(\Gamma_n)_{n\in\N}$ with the following properties:
    \begin{proplist}
        \item\label{prop:converge} There is a sequence $p_n\to p$ for which $\Gamma_n\in\mcal K_{p_n}(\mcal H)$, and
        \item\label{prop:supwhite} $\limsup_{n\to\infty}\abs{VW(\Gamma_n)}=w$, and
        \item\label{prop:distinct} The $\Gamma_n$'s are pairwise distinct.
    \end{proplist}
    Then there is a $p$-core CRG $K$ with the following properties:
    \begin{enumerate}
        \item $K\in\mcal K_p(\mcal H)$, and
        \item $g_K(p)\leq\limsup_{n\to\infty}g_{\Gamma_n}(p_n)$, and
        \item $\abs{VW(K)}=s>w$, and
        \item\label{crgseq:dalmatian} For any integer $1\leq r\leq s-w$ and any positive integer $\ell$, $\mcal F\not\mapsto K^\ell(r)$.
    \end{enumerate}
\end{theorem}

\begin{proof}
    By passing to a subsequence of $(\Gamma_n)$, we can guarantee the following additional properties:
    \begin{proplist}
        \setcounter{proplisti}{3}
        \item\label{prop:0core} $p_n\in(0,1/2)$, so $\Gamma_n\in\mcal K_0(\mcal H)$, and\hfill (\cref{prop:converge})
        \item\label{prop:wwhite} $\abs{VW(\Gamma_n)}=w$, and\hfill (\cref{prop:supwhite} and the fact that $w<\chi(\mcal F)$ is a non-negative integer)
        \item\label{prop:nblack} $\abs{VB(\Gamma_n)}\geq n$.\hfill (\cref{prop:distinct,prop:wwhite})
    \end{proplist}

    Let $\mu_n\in\Delta^{\Gamma_n}$ be the probability mass achieving $g_{\Gamma_n}(p_n)$.
    For each $n\in\N$, we can locate some integer $N_n$ such that
    \begin{proplist}[label=\roman*., ref=\roman*]
        \item\label{itm:embed} $\mu_n(x)\cdot N_n\geq n$ for all $x\in V(\Gamma_n)$, and\hfill ($\mu_n$ has full support)
        \item\label{itm:lim} $\displaystyle\absss{g_{\Gamma_n}(p_n)-{\Delta_{p_n}(\Gamma_n[\mu_n,N_n])\over\abs{\Gamma_n[\mu_n,N_n]}}}<{1\over n}$.\hfill (\Cref{prop:blowup})
    \end{proplist}
    Set $G_n\eqdef \Gamma_n[\mu_n,N_n]$.
    Note that \cref{itm:embed} implies $n\times\Gamma_n\subseteq G_n$ and that \cref{itm:lim} implies
    \[
        \limsup_{n\to\infty}{\Delta_{p_n}(G_n)\over\abs{G_n}}=\limsup_{n\to\infty}g_{\Gamma_n}(p_n).
    \]
    Furthermore,
    \begin{align*}
        \abs{\Delta_{p_n}(G_n)-\Delta_p(G_n)} &=\absss{\max_{v\in V(G_n)}d_{p_n}(v)-\max_{v\in V(G_n)}d_p(v)} \leq\max_{v\in V(G_n)}\abs{d_{p_n}(v)-d_p(v)}\leq\abs{G_n}\cdot\abs{p_n-p},
    \end{align*}
    so, since $p_n\to p$,
    \[
        \limsup_{n\to\infty}{\Delta_p(G_n)\over\abs{G_n}}=\limsup_{n\to\infty}g_{\Gamma_n}(p_n).
    \]
    \medskip

    We claim next that $K(w+1,0)\crgarr(G_n)$.
    Indeed, \cref{prop:wwhite,prop:nblack} allow us to apply \Cref{prop:manyblack} to find that $n\times K(w+1,0)\sqsubseteq n\times \Gamma_n\subseteq G_n$, where the last inclusion follows from \cref{itm:embed}.
    \medskip

    Thus, thanks to \Cref{lem:maxwhite}, we may find some $p$-core CRG $K$ such that $\abs{VW(K)}=s>w$ and $K\crgarr(G_n)$.
    Since $K\crgarr(G_n)$, we know that $\mcal F\not\mapsto K$; hence $K\in\mcal K_p(\mcal H)$.
    Furthermore,
    \[
        g_K(p)\leq\limsup_{n\to\infty}{\Delta_p(G_n)\over\abs{G_n}}=\limsup_{n\to\infty}g_{\Gamma_n}(p_n).
    \]
    \medskip

    To finish the claim, fix any $1\leq r\leq s-w$ and any $\ell\in\N$.
    In order to show that $\mcal F\not\mapsto K^\ell(r)$, it suffices to show that $K^\ell(r)\crgarr(G_n)$.

    Fix any integer $N\geq\ell$.
    Since $K\crgarr(G_n)$, there must be some $n\geq N$ for which $N\times K\sqsubseteq G_n$.
    Thus, \cref{itm:embed} allows us to apply \Cref{lem:dalmatian} to see that also $N\times K^N(r)\sqsubseteq G_n$.
    We conclude that
    \[
        N\times K^\ell(r)\subseteq N\times K^N(r)\sqsubseteq G_n,
    \]
    and so $K^\ell(r)\crgarr(G_n)$ as needed.
\end{proof}

We can now show that $\Forb(K_{t,t})$ is free of accumulation points in $(0,1]$.
\begin{proof}[Proof of \Cref{thm:ktt}]
    Set $\mcal H=\Forb(K_{t,t})$.

    To begin, we note that $p=1/2$ is never an accumulation point of any hereditary property due to~\cite[Theorem 39]{martinrandom}.

    Next, observe that if $K\in\mcal K_1(\mcal H)$, then $K$ has at most one white vertex and fewer than $\wbar\chi(K_{t,t})=t$ black vertices.
    In particular, $\mcal K_1(\mcal H)$ is finite and so $\mcal H$ has no accumulation points in $(1/2,1]$.

    Now, fix $p\in(0,1/2)$ and suppose for the sake of contradiction that $p$ is an accumulation point of $\mcal H$.
    Thus, we may find some sequence $p_n\to p$ and a sequence of CRGs $(\Gamma_n)_{n\in\N}$ for which
    \begin{itemize}
        \item $\Gamma_n\in\mcal K_{p_n}(\mcal H)$, and
        \item The $\Gamma_n$'s are distinct, and
        \item $\ed_{\mcal H}(p_n)=g_{\Gamma_n}(p_n)$.
    \end{itemize}
    Since $\ed_{\mcal H}$ is continuous, this implies that
    \[
        \lim_{n\to\infty}g_{\Gamma_n}(p_n)=\ed_{\mcal H}(p).
    \]
    Thus, we may apply \Cref{thm:crgseq} to the sequence $(\Gamma_n)$.
    Let $K$ be the $p$-core CRG guaranteed by \Cref{thm:crgseq} and note that $\abs{VW(K)}\geq 1$.
    Since $K_{t,t}$ is bipartite and $K$ has at least one white vertex, we know that $K$ has exactly one white vertex and that $\abs{VB(K)}\leq t-1$.
    Furthermore, $K\in\mcal K_p(\mcal H)$ and
    \[
        g_K(p)\leq\limsup_{n\to\infty}g_{\Gamma_n}(p_n)=\ed_{\mcal H}(p)\quad\implies\quad g_K(p)=\ed_{\mcal H}(p).
    \]
    Thus, \Cref{prop:allgray} uniquely determines $K=K(1,t-1)$.
    However, $K(1,t-1)^1(1)=K(0,t)$ and, since $\wbar\chi(K_{t,t})=t$, we know that $K_{t,t}\mapsto K(0,t)$.
    This contradicts \cref{crgseq:dalmatian} of \Cref{thm:crgseq} and so we have established the claim.
\end{proof}

Unfortunately, we have been unable to rule out the possibility of other accumulation points for $\Forb(K_{2,t})$, $t\geq 9$, though we believe that none exist.

\begin{remark}\label{rem:finite}
    Our proof of \Cref{thm:ktt} actually establishes a stronger claim.
    Set $\mcal H=\Forb(K_{t,t})$, fix $p\in(0,1]$ and fix any sequence $p_n\to p$.
    If there are CRGs $\Gamma_n\in\mcal K_{p_n}(\mcal H)$ with $\lim_{n\to\infty}g_{\Gamma_n}(p_n)=\ed_{\mcal H}(p)$, then the sequence $(\Gamma_n)$ is eventually constant.
\end{remark}

\section{Concluding remarks}\label{sec:remarks}

While we gained a better understanding of accumulation points of the edit distance function in this paper, we still have many questions. For instance, for which $p$ can a hereditary property have an accumulation point at $p$?

\begin{question}
    Let $\mcal A$ be the set of all $p\in[0,1/2]$ for which there is some non-trivial hereditary property with an accumulation point at $p$.
    What can be said about the set $\mcal A$?
    In particular:
    \begin{itemize}
        \item What is $\sup\mcal A$?
        \item Is $\mcal A$ an interval?
    \end{itemize}
\end{question}
Currently, all that we know about the set $\mcal A$ is that $0,1/4\in\mcal A$ and that $\sup\mcal A\leq 1-\varphi^{-1}=0.382\ldots$ where $\varphi=1.618\ldots$ is the golden ratio.
The latter follows from the work of Martin and Riasanovsky~\cite[Theorem 39]{martinrandom}.
The work of Martin--Riasanovsky additionally suggests that perhaps $\sup\mcal A\leq 1/3$ (c.f.\ \cite[Section 2.3 \& Proposition 35]{martinrandom}).

We expect that, at the very least, $1/n\in\mcal A$ for all integers $n\geq 4$.
Indeed, one could likely build on the ideas in the proof of \Cref{thm:otheracc} and construct a family $\mcal F$ of $(n-2)$-regular graphs so that the property $\Forb(\{K_{1,n}\}\cup\mcal F)$ has an accumulation point at $p=1/n$.
\medskip

Beyond this, it is worth pointing out that our construction of a hereditary property with an accumulation point at $p=1/4$ from \Cref{thm:otheracc} required an infinite set of forbidden graphs.
It is natural to wonder if this was necessary.
\begin{question}
    Suppose that $\mcal H=\Forb(\mcal F)$ where $\mcal F$ is finite.
    Can $\mcal H$ have any accumulation points in $(0,1)$?
\end{question}

We expect a negative answer, though, since this is the first systematic investigation into the accumulation points of the edit distance function, our current knowledge is severely limited.
Regardless of the answer to the above question, we are confident to conjecture the following extension of \Cref{thm:ktt}.
\begin{conj}
    For any $s,t\in\N$, $\Forb(K_{s,t})$ has no accumulation points in the interval $(0,1]$.
\end{conj}
We note that it is necessary to go beyond the ideas used in \Cref{thm:ktt} in order to establish this conjecture.
Indeed, for every $k\geq 4$, Marchant--Thomason~\cite[Theorem 3.27]{MarchantThomasonExtremalColors} constructed infinitely many $1/k$-core CRGs achieving $\ed_{\Forb(K_{1,k})}(1/k)$; therefore, the analogue of \Cref{rem:finite} fails to hold in these cases.

Despite this, the same reasoning used in \Cref{thm:ktt} implies that $\Forb(K_{s,t})$ is free of accumulation points in $[1/2,1]$ and that if $p\in(0,1/2)$ happens to be an accumulation point, then $K(1,\min\{s,t\}-1)$ achieves $\ed_{\Forb(K_{s,t})}(p)$.
\medskip

We furthermore believe that accumulation points are relatively rare.
\begin{conj}
    Any non-trivial hereditary property has only finitely many accumulation points.
\end{conj}

Finally, we believe that accumulation points can occur only in one half of the interval $(0,1)$.
\begin{conj}
    For any non-trivial hereditary property, either $(0,1/2]$ or $[1/2,1)$ is free of accumulation points.
\end{conj}

\bibliographystyle{abbrv}
\bibliography{references}
\end{document}